\input amstex
\documentstyle{amsppt}
\loadeusm
\topmatter
\title
Frobenius Pseudoprimes
\endtitle
\date March 1999
\enddate
\email grantham\@super.org
\endemail
\author Jon Grantham \endauthor
\endsubjclass
\address Institute for Defense Analyses, Center for Computing Sciences,
17100 Science Drive, Bowie, MD 20715 
\endaddress
\define\mepseudotwo{11}
\define\mepseudothree{12}
\abstract
The proliferation of probable prime tests in recent years has
produced a plethora of definitions with the word
``pseudoprime'' in them.  Examples include pseudoprimes, Euler
pseudoprimes, strong pseudoprimes, Lucas
pseudoprimes, strong Lucas pseudoprimes, extra strong Lucas
pseudoprimes and Perrin pseudoprimes.  Though these tests represent a
wealth of ideas, they exist as a hodge-podge of definitions rather
than as examples of a more general theory.
It is the goal of this paper to
present a way of viewing many
of these tests as special cases of a general principle, as
well as to re-formulate them in the context of finite fields.

One aim of the reformulation is to enable the creations of
stronger tests;
another is to aid in proving results about large classes
of pseudoprimes.
\endabstract

\NoBlackBoxes
\def\as{2}
\def\adams{1}
\def\agp{3}
\def\agptwo{4}
\def\arno{5}
\def\atkin{6}
\def\bw{7}
\def\gord{8}
\def\gordpom{9}
\def\mediss{10}
\def\gurak{13}
\def\guy{14}
\def\jac{15}
\def\ksw{16}
\def\len{17}
\def\mojones{18}
\def\mon{19}
\def\dopo{20}
\def\pom{21}
\def\psw{22}
\def\rabin{23}
\def\robin{24}
\def\rot{25}
\def\rots{26}
\def\sze{27}
\def\gcmd{\operatorname{gcmd}}
\def\disc{\operatorname{disc}}
\endtopmatter
\document
\head{\S 1 Introduction}
\endhead
Fermat's Little Theorem tells us that
$a^{p-1}\equiv 1\bmod p$
for $p$ an odd prime; thus we have an easy way to prove that many
numbers are composite.  For example, since $2^{90}\equiv 64$ mod $91$,
we prove that $91$ is composite.  The technique of
repeated squaring can be used to perform the required exponentiation very
rapidly.

This test is not foolproof.  In particular,
$2^{340}\equiv 1\bmod 341$.  Composites which fool the test with
$a=2$ are called pseudoprimes, and in general, composites $n$ with
$a^{n-1}\equiv 1\bmod n$ are pseudoprimes to the base $a$.  

The existence of such numbers provides incentive to create other tests
which are similarly fast, but which may have fewer, or at least
different, ``pseudoprimes.''  Two of these tests are
more elaborate versions of the test described above and create the notions of
Euler pseudoprime and strong pseudoprime.

Most other tests, however,  involve recurrence sequences.  
One reason that pseudoprimes based on recurrence sequences have
attracted interest is that the pseudoprimes for these sequence are
often different from ordinary pseudoprimes.  In fact, nobody has
claimed the \$620 offered for a Lucas pseudoprime with parameters
$(1,-1)$ (see \S 2 for a definition of this term), congruent to $2$
or $3\bmod 5$, that is also a pseudoprime to the base $2$ \cite{\psw, \guy}.

Furthermore, some tests based on higher order recurrence sequences seem to
have few pseudoprimes.  Adams and Shanks \cite{\as} introduced such
a test based on a third order recurrence sequence known as Perrin's
sequence.

A problem with tests based on recurrence sequences is that
analysis of the tests can be difficult.  For example, the concepts
of Lucas and Lehmer pseudoprimes have been analyzed separately in
the literature.  In \S 2, we show that they are
equivalent definitions.

Instead of recurrence sequences, the language used in this paper
is that of finite fields.  In particular, when $n$ is a prime,
$f(x)\in\Bbb{Z}[x]$,  and $n$ does not divide $\disc(f)$, the residue ring
$\Bbb{Z}[x]/(n,f(x))$ is  a  product of finite
fields of characteristic $n$.  When $n$ is composite, this ring is
not equal to such a product.
 For a given composite $n$, this fact is often easy to
discover, thus providing a quick proof of the compositeness of $n$.
Using properties of finite fields, we establish the definition of
Frobenius probable prime, which is a generalization, and
sometimes a strengthening, of many existing definitions.

We introduce the concept of Frobenius pseudoprimes, not to inflict a new
and different notion of pseudoprimality on the mathematical world, but
to show that many existing pseudoprimality tests can be generalized and
described in terms of finite fields.  In fact, we show that some specific
instances of the Frobenius pseudoprime test are equivalent to 
other pseudoprimality tests, and stronger than many of them.

In \cite{\mepseudotwo}, we use the structure given by the introduction of finite
fields to show that the probability of error in declaring a number
``prime'' using a certain
Frobenius test is less than $\frac 1{7710}$.  In 1980, Monier
\cite{\mon} and Rabin \cite{\rabin} proved that the Strong Probable
Prime Test has probability of error at most $\frac 14$.
Although the test introduced in \cite{\mepseudotwo}
 has
asymptotic running time three times that of the Strong Probable Prime
Test, the proven bound on the error is much smaller than the $\frac 1
{64}$ achieved through three Strong Probable Prime Tests.

Perhaps the primary benefit of this approach is that instead of having
to prove ten different theorems about ten different types of
pseudoprimes, one can prove one theorem about Frobenius pseudoprimes
and apply it to each type of pseudoprime.  

In \cite{\mepseudothree}, the techniques of \cite{\agp} are used to prove that
for any monic, squarefree polynomial,
there are infinitely many Frobenius pseudoprimes.   In particular,
this theorem answers a 1982 conjecture of Adams and Shanks \cite{\as} that
there are infinitely many Perrin pseudoprimes.  It also proves the
infinitude of the types of pseudoprimes defined by Gurak \cite{\gurak}
and Szekeres \cite{\sze}.

We should note that the idea of primality testing in finite fields is
not entirely new.  Lenstra's Galois Theory Test \cite{\len} is a
method of proving primality using finite fields.  In \cite{\mediss}, I
describe the relation between the two ideas.  The combination of
finite fields and pseudoprimes also exists implicitly in some other
works, such as \cite{\sze}.  The goal here, however, is different.  I
am trying to provide a clear theoretical framework in which various
existing probable prime tests can be generalized and analyzed.

\head{\S 2 A Wealth of Pseudoprimes}
\endhead
For the purposes of this paper, 
the following test for primality will be considered foolproof.
If an integer is denoted by the letter $p$, then $p$ is prime.
If $q$ is a prime power, we let $\Bbb{F}_q$ denote a finite
field with $q$ elements.

We begin by reviewing many of the existing notions of
pseudoprimality.  Each of these definitions of ``pseudoprime''
characterizes composite numbers with a certain property.  In each of
these cases, it can be proven that all prime numbers (with a finite, known
set of exceptions) have this property.

This paper does not pretend to be an exhaustive treatment of all
notions of pseudoprimality.  For example, nothing is said about
elliptic pseudoprimes \cite{\gord}.  

Fermat's Little Theorem tells us that if $p$ is prime, then
$a^{p-1}\equiv 1 \bmod p$, if $p\nmid a$.
The original notion of a pseudoprime (sometimes called a Fermat
pseudoprime) involves counterexamples to the converse of this theorem.

\proclaim{Definition}
A {\bf pseudoprime} to the base $a$ is a composite number $n$ such
$a^{n-1}\equiv 1\bmod n$.  
\endproclaim

\proclaim{Definition}
A number $n$ which is a pseudoprime to all bases $a$ with $(a,n)=1$ is
a {\bf Carmichael number}.
\endproclaim
We also know that if $p$ is an odd prime, then 
$\topsmash{a^{(p-1)/2}\equiv\left(a\over p\right)} \bmod p$,  where 
$\topsmash{\left(a\over
p\right)}$ is the Jacobi symbol.  The converse of this theorem leads to
the definition of Euler pseudoprime, due to Raphael Robinson
\cite{\robin}.

\proclaim{Definition}
An {\bf Euler pseudoprime} to the base $a$ is an odd composite number $n$ with 
$(a,n)=1$
such that $a^{(n-1)/2}\equiv\left(a\over n\right) \bmod n$.
\endproclaim

An Euler pseudoprime to the base $a$ is also a pseudoprime to the base
$a$.

If $n\equiv 1 \bmod 4$, we can also look at $a^{(n-1)/{2^k}}$ for $k>1$,
and doing so gives us the definition of strong pseudoprime
\cite{\psw}, due independently to R. Dubois and John Selfridge.

\proclaim{Definition}
A {\bf strong pseudoprime} to the base $a$ is an odd composite $n=2^rs+1$
with $s$ odd such that either $a^s\equiv 1 \bmod n$, or $a^{2^ts}\equiv -1$
for some integer $t$, with $r>t\ge 0$.
\endproclaim

A strong pseudoprime to the base $a$ is also an Euler pseudoprime to
the base $a$  \cite{\mon, \psw}.

It is possible to define notions of pseudoprimality based on
congruence properties of recurrence sequences.  The simplest of these
are based on the Lucas sequences $U_n(P,Q)$, where $P$ and $Q$
are integers, $U_0=0$,
$U_1=1$ and $U_n=PU_{n-1}-QU_{n-2}$.
(When $P=1$ and $Q=-1$, this is the Fibonacci sequence.)

We recall the fact that we can express $U_n$ in terms of roots of the
polynomial $f(x)=x^2-Px+Q$.  If $\alpha$ and
$\beta$ are roots of $f(x)$ in a commutative ring (with identity),
 with $\alpha-\beta$ invertible,
then $U_n=(\alpha^n-\beta^n)/(\alpha-\beta)$.
By induction on $n$, this equality holds
even if there are more than two distinct roots of $f(x)$.

\proclaim{Theorem 2.1}
Let $U_n=U_n(P,Q)$ and $\Delta=P^2-4Q$.  If $p\nmid 2Q\Delta$, then
$\botsmash{U_{p-\left(\Delta\over p\right)}}\equiv 0 \bmod p$.
\endproclaim
\demo{Proof}
Since $p\nmid\Delta$, $x^2-Px+Q$ has distinct roots, $\alpha$ and
$\beta$, in $\bar\Bbb{F}_p$.

If $\left(\Delta\over p\right)=1$, then $f(x)$ factors mod $p$, and
$\alpha$ and $\beta$ are in $\Bbb{F}_p$.  Thus $\alpha^{p-1}\equiv
\beta^{p-1}\equiv 1\bmod p$.  So $U_{p-1}\equiv (1-1)/(\alpha-\beta)=0\bmod
p$.

If $\left(\Delta\over p\right)=-1$, then $f(x)$ does not factor, and
the roots of $f(x)$ lie in $\Bbb{F}_{p^2}$.  The
Frobenius automorphism permutes the roots of $f(x)$, so
$\alpha^p\equiv\beta$ and $\beta^p\equiv\alpha$.  Thus
$U_{p+1}\equiv(\alpha\beta-\beta\alpha)/(\alpha-\beta)=0\bmod p$.
\enddemo

\proclaim{Definition}
Let $U_n=U_n(P,Q)$ and $\Delta=P^2-4Q$.
A {\bf Lucas pseudoprime} with parameters $(P,Q)$ is a composite $n$
with $(n,2Q\Delta)=1$ such that $U_{n-\left(\Delta\over n\right)}\equiv 0 \bmod
n$. 
\endproclaim

Baillie and Wagstaff \cite{\bw}
gave a version of this test that is analogous to the
strong pseudoprime test.  We first define the sequence $V_n(P,Q)$ to
be the sequence with $V_0=2$, $V_1=b$, and $V_n=PV_{n-1}-QV_{n-2}$.
Note that $V_n=\alpha^n+\beta^n$, where $\alpha$ and $\beta$ are
distinct roots of $x^2-Px+Q$.

\proclaim{Theorem 2.2}
Let $U_n=U_n(P,Q)$ and $\Delta=P^2-4Q$.  Let $p$ be a prime not
dividing $2Q\Delta$.  Write  $p=2^rs+\left(\Delta\over p\right)$, where
$s$ is odd.   Then either $U_s\equiv 0$ or $V_{2^ts}\equiv 0 \bmod p$
for some $t$, $0\le t < r$.
\endproclaim
\demo{Proof}
As above, since $p\nmid\Delta$, $x^2-Px+Q$ has distinct roots
in $\bar\Bbb{F}_p$.
If $\left(\Delta\over p\right)=1$, the roots $\alpha$ and $\beta$ are
elements of $\Bbb{F}_p$.  If $\botsmash{\left(\Delta\over p\right)}=-1$,
they are elements of $\Bbb{F}_{p^2}$.  
Their product is $\alpha\beta=Q\neq 0$, so $\beta\neq 0$.

We have that
$U_{2^rs}=(\alpha^{2^rs}-\beta^{2^rs})/(\alpha-\beta)\equiv 0 \bmod p$.  Thus
$(\alpha/\beta)^{2^rs}\equiv 1$  mod $p$.  Taking square
roots, we have that either $(\alpha/\beta)^{2^ts}\equiv -1$
for some $t$ with $0\le t<r$, or $(\alpha/\beta)^s\equiv 1$. 

In the first case, we can conclude $V_{2^ts}\equiv 0 \bmod p$, while in
the second case we have that $U_s\equiv 0 \bmod p$.
\enddemo

\proclaim{Definition}
Let $U_n=U_n(P,Q)$, $V_n=V_n(P,Q)$, and $\Delta=P^2-4Q$.
A {\bf strong Lucas pseudoprime} with parameters $(P,Q)$ is a
composite $n=2^rs+\left(\Delta\over n\right)$, where $s$ is odd and
$(n,2Q\Delta)=1$ such that either $U_s\equiv 0 \bmod n$ or $V_{2^ts}\equiv 0
\bmod n$ for some $t$, $0\le t<r$.
\endproclaim

Any strong Lucas pseudoprime is also a Lucas pseudoprime with the same
parameters.

Jones and Mo \cite{\mojones} have recently given another test that
relies on the sequences $U_n(b,1)$ and $V_n(b,1)$.

\proclaim{Theorem 2.3}
Let $U_n=U_n(b,1)$, $V_n=V_n(b,1)$,
 and $\Delta=b^2-4$.  Let $p$ be a prime not
dividing $2\Delta$.  Write  $p=2^rs+\left(\Delta\over p\right)$, where
$s$ is odd.   Then either $U_s\equiv 0 \bmod p$ and $V_s\equiv\pm 2
\bmod p$, or $V_{2^ts}\equiv 0 \bmod p$,
for some $t$, $0\le t < r-1$.
\endproclaim
\demo{Proof}
By Theorem 2.2, it suffices to show that $V_{2^{r-1}s}\not\equiv 0$,
and that if $U_s\equiv 0$, then $V_s\equiv\pm 2$.

Note that $V_n=\alpha^n+\alpha^{-n}$, where $\alpha$ is a root of
$x^2-bx+1$.  So $V_{2^{r-1}s}
\equiv 0 \bmod p$ implies $\alpha^{2^rs}\equiv -1$.  
If 
$\left(\Delta\over p\right)=1$, then $\alpha\in\Bbb{F}_p$, and we have
a contradiction. If
$\left(\Delta\over p\right)=-1$, then
$\alpha^p\equiv\alpha^{-1}$. So $\alpha^{p+1}
\equiv 1\not\equiv -1$, and we also have a contradiction.

If $U_s\equiv 0 \bmod p$, then $\alpha^s\equiv\alpha^{-s} \bmod p$, and thus
$\alpha^{2s}\equiv 1$.  We must have $\alpha^s\equiv\pm 1$, and thus
$V_s\equiv\pm 2$.
\enddemo

\proclaim{Definition}
Let $U_n=U_n(b,1)$, $V_n=V_n(b,1)$, and $\Delta=b^2-4$.
An {\bf extra strong Lucas pseudoprime} to the base $b$ is a composite
$n=2^rs+\left(\Delta\over n\right)$, where $s$ is odd and $(n,2\Delta)=1$,
such that either $U_s\equiv 0 \bmod n$ and $V_s\equiv\pm 2 \bmod n$,
or $V_{2^ts}\equiv 0 \bmod n$ for some $t$ with $0\le t<r-1$.
\endproclaim

Any extra strong Lucas pseudoprime base $b$ is a strong Lucas pseudoprime
with parameters $(b,1)$.

The definition of Lehmer pseudoprime, due to Rotkiewicz \cite{\rot},
is related to that of Lucas pseudoprime.

The Lehmer sequence with parameters $(L,Q)$ is defined by $\bar U_0=0$,
$\bar U_1=1$, $\bar U_k=L\bar U_{k-1}-Q\bar U_{k-2}$ for $k$ odd, and $\bar
U_k=\bar U_{k-1}-Q\bar U_{k-2}$
for $k$ even.
It is easy to see by induction that
$\bar U_k=(\alpha^k-\beta^k)/(\alpha^2-\beta^2)$ for $k$ even, and $\bar U_k=
(\alpha^k-\beta^k)/(\alpha-\beta)$ for $k$ odd, where $\alpha$ and 
$\beta$ are roots of $x^2-\sqrt{L}x+Q$.

\proclaim{Definition}
Let $D=L-4Q$ and $\epsilon(n)=\left(LD\over n\right)$.
A {\bf Lehmer pseudoprime} with parameters $(L,Q)$ is
a composite $n$ with $(2LD,n)=1$ and $\bar U_{n-\epsilon(n)}\equiv 0$  mod
$n$.
\endproclaim

Lehmer pseudoprimes can be analyzed by the same means as Lucas
pseudoprimes, because of the following new result.

\proclaim{Theorem 2.4}
An integer $n$ is a Lehmer pseudoprime with parameters $(L,Q)$ if and
only if it is a Lucas pseudoprime with parameters $(L,LQ)$.
\endproclaim
\demo{Proof}
Let $D=L-4Q$ be the discriminant of the Lehmer sequence.
The characteristic polynomial of the Lucas sequence is $f(x)=x^2-Lx+LQ$,
which has discriminant $L^2-4LQ=LD$.  So
$\epsilon(n)=\left(L^2-4LQ\over n\right)$.  The roots of $f(x)$ are
$\alpha=\frac{L+\sqrt{L^2-4LQ}}2$ and
$\beta=\frac{L-\sqrt{L^2-4LQ}}2$. 
The characteristic polynomial of the Lehmer sequence is
$g(x)=x^2-\sqrt{L}x+Q$.  Its roots are $\alpha'=\alpha/\sqrt{L}$ and
$\beta'=\beta/\sqrt{L}$.

Thus $L^{(n-\epsilon(n))/2}\bar U_{n-\epsilon(n)}=U_{n-\epsilon(n)}$,
and we conclude that $\bar U_{n-\epsilon(n)}\equiv 0\bmod n$ if and only if
$U_{n-\epsilon(n)}\equiv 0\bmod n$.  This proves the theorem.

\enddemo

Rotkiewicz \cite{\rots}
has also given a definition of strong Lehmer pseudoprime.

\proclaim{Definition}
Let $\bar U_k$ be as in the definition of Lehmer pseudoprime.  Let $\bar
V_n$ satisfy $\bar V_0=2$, $\bar V_1=1$, $\bar V_k=L\bar V_{k-1}-Q\bar
V_{k-2}$ for $k$ even, and 
$\bar V_k=\bar V_{k-1}-Q\bar V_{k-2}$ for $k$ odd.  
Let $\epsilon(n)$ be as above.
An odd composite number $n=2^rs+\epsilon(n)$ 
is a {\bf strong Lehmer pseudoprime} with parameters $(L,Q)$ if $(n,DQ)=1$
and either $\bar U_s\equiv 0 \bmod n$ or $\bar V_{2^ts}\equiv 0$ for some $t$
with $0\le t<r$.
\endproclaim

\proclaim{Theorem 2.5}
An integer $n$ is a strong Lehmer pseudoprime with parameters $(L,Q)$ if and
only if it is a strong Lucas pseudoprime with parameters $(L,LQ)$.
\endproclaim
\demo{Proof}
The technique of proof is exactly the same as for Theorem 2.4.
\enddemo

In a series of papers in the 1980s (\cite{\as}, \cite{\ksw} and \cite{\adams}), 
Adams, Shanks and co-authors proposed and analyzed a
pseudoprime test based on a third order recurrence sequence known as
Perrin's sequence, and on generalizations of this sequence.
A good exposition of the test is given in \cite{\arno}.

We consider sequences $A_n=A_n(r,s)$ defined by the following relations:
$A_{-1}=s$, $A_0=3$, $A_1=r$, and $A_n=rA_{n-1}-sA_{n-2}+A_{n-3}$.
Let $f(x)=x^3-rx^2+sx-1$ be the associated polynomial and $\Delta$ its 
discriminant.
Perrin's sequence is $A_n(0,-1)$.

\proclaim{Definition}
The {\bf signature mod $m$} of an integer $n$ with respect to the sequence
$A_k(r,s)$ is the $6$-tuple
$(A_{-n-1},A_{-n},A_{-n+1},A_{n-1},A_n,A_{n+1}) \bmod m$.
\endproclaim

\proclaim{Definitions}
An integer $n$ is said to have an {\bf S-signature} if its signature  mod $n$
is congruent to $(A_{-2},A_{-1},A_0,A_0,A_1,A_2)$.

An integer $n$ is said to have a {\bf Q-signature} if its signature
  mod $n$ is congruent to $(A,s,B,B,r,C)$, where for some integer $a$
with $f(a)\equiv 0 \bmod n$, $A\equiv a^{-2}+2a$, $B\equiv
-ra^2+(r^2-s)a$, and $C\equiv a^2+2a^{-1}$.

An integer $n$ is said to have an {\bf I-signature} if its signature
mod $n$ is congruent to $(r,s,D',D,r,s)$, where $D'+D\equiv rs-3$ mod
$n$ and $(D'-D)^2\equiv \Delta$.
\endproclaim

\proclaim{Definition}
A {\bf Perrin pseudoprime} with parameters $(r,s)$ is an odd composite $n$
and either 

1) $\left(\Delta\over n \right)=1$ and $n$ has an S-signature or an 
I-signature, or

2) $\left(\Delta\over n \right)=-1$ and $n$ has a Q-signature.
\endproclaim

In the past, the term Perrin pseudoprime has referred only to
pseudoprimes with respect to Perrin's original sequence, but we feel it
is useful to have a convenient name for composites having an
acceptable signature for other such sequences.

Also, we have omitted a portion of the test involving quadratic forms.
If $n$ has an I-signature, it is possible to 
construct a quadratic form representing $n$.
A prime with an I-signature can only be represented by forms lying in
a certain subgroup of the class group of quadratic forms with
discriminant $\Delta$.  The only examples found where
this portion of the test has exposed composites involve
pseudoprimes very small compared to the discriminant of the
associated polynomials.  In fact, the
polynomials were cleverly constructed specifically to have
pseudoprimes which would be exposed in this way.
In particular, 
there are no known examples where the quadratic form exposes a
composite for the test using Perrin's sequence.
The test does not apply to integers with Q-signatures or
S-signatures.
The interested reader should consult \cite{\as} for details.

Generalizations to higher order recurrence sequences have been given
by Gurak \cite{\gurak}.  His basic definition is shown below to be
subsumed in the definition of Frobenius pseudoprimes.  Later in his
paper, he gives ideas as to how his test could be made stronger.  He
does not, however, give exact definitions of other notions of
pseudoprimality.  

A nice variation on these tests is given by Szekeres \cite{\sze}.

\proclaim{Definition}
Let $f(x)$ be an irreducible polynomial in $\Bbb{Z}[x]$ and
$\beta_1,\dots,\beta_k$ be its roots.
A {\bf pseudoprime (in the sense of Szekeres)} is a composite $n$ such that
for every symmetric polynomial $S(x_1,\dots,x_k)\in\Bbb{Z}[x_1,\dots,x_k]$,
$S(\beta_1^n,\dots,\beta_k^n)\equiv S(\beta_1,\dots,\beta_k)$  mod $n$.
\endproclaim

Szekeres does not test signatures and notes, ``Signatures can be
tested without much additional effort but they don't seem to add
significantly to the efficiency of primality testing through higher
order Lucas sequences\dots''  He also does not use knowledge gained
from Jacobi symbols in his test.

Atkin has proposed a specific test based on arithmetic modulo
polynomials; it shares some similarities with the test described in
Section 3.  He describes it fully in \cite{\atkin}.

\head{\S 3 Frobenius Pseudoprimes}
\endhead
In this section, We will be introducing the definition of Frobenius
pseudoprime.  This definition does not in and of itself solve any open
questions in the subject.  We do, however, aim to provide a clearer way
of thinking about the definitions given in the previous section.  Some
open questions are solved in \cite{\mepseudotwo} and
\cite{\mepseudothree}.

We first prove some elementary facts about polynomials over finite fields.
In particular, we exploit the following fact.  Given a polynomial
$f(x)$ of degree $d$, we can factor it as $\prod_{1}^d F_i(x)$,
where each $F_i(x)$ is the product of the irreducible polynomials of degree $i$
dividing $f(x)$.

More precisely, we define these polynomials as follows.
Let $f_0(x)=f(x)$.
For $1\le i\le d$, define $F_i(x)=\gcd(x^{p^i}-x,f_{i-1}(x))$ in
$\Bbb{F}_p[x]$ and
$f_i(x)=f_{i-1}(x)/F_i(x)$.

\proclaim{Theorem 3.1}
Let $p$ be an odd prime, and let $f(x)$ be a monic polynomial in
$\Bbb{F}_p[x]$ of degree $d$ with discriminant $\Delta$.
Assume $p\nmid f(0)\Delta$.

1)  
WE have $f_d(x)=1$, and for each $i$,  $1\le i\le d$, $i|\deg(F_i)$.

2) For $2\le i\le d$, $F_i(x)|F_i(x^{p})$.

3)  Let $S=\sum\Sb{2|i}\endSb\deg(F_i(x))/{i}$.
Then $(-1)^S=\left(\Delta\over p\right)$.   That is, 
if $\left(\Delta\over p\right)=1$, then $S$ is an even integer, and
if $\left(\Delta\over p\right)=-1$, then $S$ is an odd integer.
\endproclaim
\demo{Proof}
1) The polynomial
$x^{p^i}-x$ is the product of all irreducible polynomials in $\Bbb{F}_p[x]$
with degree dividing $i$.   Inductively, we see that
$$F_i(x)=\gcd(\prod_{j|i}(x^{p^j}-x)^{\mu(i/j)},f(x)),$$
where $\mu$ is the Mobius function.
 Thus $F_i(x)$ is the product of all irreducible
factors of $f(x)$ of degree exactly $i$, and $f_i(x)$ is the
product of the irreducible polynomials dividing
$f(x)$ with degree greater than $i$.

Since $\Delta\not=0$, $f(x)$ is squarefree, and
since $f(x)$ is equal to the product of its irreducible
factors, $f(x)=\prod_{1\le
i\le d} F_i(x)$, so $f_d(x)=f(x)/(\prod_{1\le i\le d} F_i(x))=1$.

2) In fact, for any nonzero polynomial $g(x)\in\Bbb{F}_p[x]$, we have
$g(x)|g(x^p)$ since $g(x^p)=g(x)^p$.

3) The degree of $F_i(x)$ is $i$ times the number of irreducible
factors of $f$ of degree $i$.
So $S$ is equal to the number of irreducible factors of $f$ of even degree.

If $f(x)$ is irreducible mod $p$, then
$d=\deg(f)$ is the least power with $\alpha^{p^d}=\alpha$, and the map
$\alpha\mapsto\alpha^p$ has order equal to $d$.  Thus, that map is
a generator of the Galois group of $f$ over $\Bbb{F}_p$.  

For all polynomials, the Galois group acts
transitively on the roots of each irreducible factor of $f$ over
$\Bbb{F}_p$.  Thus, $S$ gives
the number of cycles of even length in the Frobenius automorphism.
Cycles of even length are odd (and vice versa), so the parity of $S$
determines whether the automorphism is odd or even.  Since the
discriminant is the product of the square of the
differences of the roots of $f(x)$,
this parity is
also determined by $\left(\Delta\over p\right)$.    For a more
detailed proof of this fact, see \cite{\jac}.
\enddemo

As an example, let $f(x)=x^4+12x+1$.  (It is irreducible over $\Bbb
Q$.)  Let $p=89$.
We have $x^{89}-x\equiv 59x^3+51x^2+20x+86 \bmod (89,f(x))$, so
$$F_1(x)=\gcd(f(x), 59x^3+51x^2+20x+86)=x+78,$$ and
$f_1(x)=x^3+11x^2+32x+8$.

Since $x^{89^2}-x\equiv 64x^2+86x+19 \bmod (89,f_1(x))$, and
$$F_2(x)=\gcd(f_1(x), 64x^2+86x+19)=1,$$ 
we have $f_2(x)=f_1(x)=x^3+11x^2+32x+8$.

Next, $x^{89^3}-x\equiv 0 \bmod (89,f_2(x))$, so $F_3(x)=f_2(x)$ and $f_3(x)=1$.

Thus $F_4(x)=f_4(x)=1$.  

Note that $x^{89}\equiv 25x^2+x+59 \bmod (89,F_3(x))$.  We verify that $F_3(25x^2+x+59)\equiv 0 \bmod (89,F_3(x))$.

Finally, the discriminant of $f(x)$ is $-559616$.  $\left({-559616}\over
{89}\right)=1$, which agrees with the fact that $S=0$.

\vskip.1in

We would like to define any composite with satisfies the consequences
of this theorem to be a type of pseudoprime, but we may run into a problem
when we take the gcd of two polynomials modulo a composite, since we
are working over a ring that is not a domain.  

With this in mind, we consider the following definition.

\proclaim{Definition}
Let $f(x),g_1(x),g_2(x)$ be monic polynomials over a commutative ring
(with identity).
We say that polynomial $f(x)$ is the {\bf greatest common monic
divisor (gcmd)} of $g_1(x)$ and $g_2(x)$ if the ideal generated by $g_1(x)$
and $g_2(x)$ is equal to the ideal generated by $f(x)$.  We write
$f(x)=\gcmd(g_1(x),g_2(x))$. (Note that $\gcmd(g_1(x),g_2(x))$ does not
necessarily exist.)
\endproclaim

Clearly if $\gcmd(g_1(x),g_2(x))$ exists, it is a common monic divisor
of $g_1(x)$ and $g_2(x)$ of greatest degree.  Further, it is not hard
to show that the $\gcmd$ is unique.  The $\gcmd$ has the following additional
property.

\proclaim{Proposition 3.2} Let $p|n$ and $g_1(x),g_2(x)$ be monic
polynomials in $\Bbb{Z}[x]$.
If $f(x)=\gcmd(g_1(x),g_2(x))$ in $(\Bbb{Z}/n\Bbb{Z}[x])$, then $f(x)\equiv
\gcd(g_1(x),g_2(x)) \bmod p$, where the $\gcd$ is taken in
$(\Bbb{Z}/p\Bbb{Z})[x]$.
\endproclaim
\demo{Proof}
For $i=1,2$, we have that $g_i(x)\equiv k_i(x)f(x) \bmod n$, 
for some $k_i(x)\in\Bbb{Z}[x]$.  Thus $f(x)|g_i(x)$ in $(\Bbb{Z}/p\Bbb{Z})[x]$.

We have that $f(x)\equiv g_1(x)h_1(x)+g_2(x)h_2(x) \bmod n$ for some
$h_1(x),h_2(x)\in\Bbb{Z}[x]$.
Thus $f(x)\equiv g_1(x)h_1(x)+g_2(x)h_2(x) \bmod p$, and by the
definition of $\gcd$, $f(x)=\gcd(g_1(x),g_2(x))$ in
$(\Bbb{Z}/p\Bbb{Z})[x]$.
\enddemo

\proclaim{Corollary 3.3}
If $\gcmd(g_1(x),g_2(x))$ exists in $(\Bbb{Z}/n\Bbb{Z})[x]$, then
for all $p$ dividing $n$,
$\gcd(g_1(x),g_2(x))$ has the same degree.
\endproclaim
\demo{Proof}
Since the leading coefficient of the $\gcmd$ is $1$, that coefficient
is the leading coefficient of all the $\gcd$s, by Proposition 3.2.
Thus, they all have the same degree.
\enddemo

\proclaim{Proposition 3.4}
Assume, for each $p|n$, $\gcd(f(x),g(x))=1$ in $\Bbb{F}_p[x]$.
Then $\gcmd(f(x),g(x))=1$ in $(\Bbb{Z}/n\Bbb{Z})[x]$.
\endproclaim
\demo{Proof}
By the Chinese Remainder Theorem, it suffices to prove that 
the $\gcmd$ of $f(x)$ and $g(x)$ is $1$
in $(\Bbb{Z}/p^j\Bbb{Z})[x]$ for any integer $j\ge 1$.
We proceed by induction.  We know that there exist polynomials
$a_1(x)$, $b_1(x)$, and $k_1(x)$ such that 
$f(x)a_1(x)+g(x)b_1(x)=1+pk_1(x)$.  Assume that there exist
polynomials $a_j(x)$, $b_j(x)$ and $k_j(x)$ with
$f(x)a_j(x)+g(x)b_j(x)=1+p^jk_j(x)$.
Multiplying by
$pk_1(x)$,
$pf(x)a_j(x)k_1(x)+pg(x)b_j(x)k_1(x)=pk_1(x)+p^{j+1}k_1(x)k_j(x)$.

Substituting for $pk_1(x)$,
$pf(x)a_j(x)k_1(x)+pg(x)b_j(x)k_1(x)
=-1+f(x)a_1(x)+g(x)b_2(x)+p^{j+1}k_1(x)k_j(x)$.
Rearranging,
$f(x)[a_1(x)-pa_j(x)k_1(x)]+g(x)[b_1(x)-pb_j(x)k_1(x)]
=1+p^{j+1}[-k_1(x)k_j(x)]$.

Thus we have shown that $\gcmd(f(x),g(x))=1$ in
$(\Bbb{Z}/p^{j+1}\Bbb{Z})[x]$, and the Proposition is proven.
\enddemo

The concept of $\gcmd$ would not be useful in the context of this
paper if it were
difficult to calculate.  We have the following result, which aids us
in testing primality.

\proclaim{Proposition 3.5}
The Euclidean algorithm will either find the $\gcmd$ of two monic
polynomials in $(\Bbb{Z}/n\Bbb{Z})[x]$ or find a proper factor of $n$.
\endproclaim
\demo{Proof}
The Euclidean algorithm will only fail to finish if one of the
divisions fails due to the leading coefficient
of a non-zero remainder being non-invertible.  
This coefficient will have a $\gcd$ with $n$ that is a nontrivial factor of $n$.

If the Euclidean algorithm terminates (i.e., one of the remainders is
zero), we have inductively that the last non-zero remainder is a
divisor of the two polynomials and can be written as a linear
combination of the two.  The proof is the same as the proof of
correctness of the Euclidean
algorithm over $\Bbb{F}_p[x]$.  Since the leading coefficient of the
last non-zero remainder is invertible, this remainder
can be made monic by division, and we find the $\gcmd$.
\enddemo

\proclaim{Definition}
Let $f(x)\in\Bbb{Z}[x]$ be a monic polynomial of degree $d$
with discriminant $\Delta$.   
An odd integer $n>1$ is said to pass the {\bf Frobenius
probable prime test} with respect to $f(x)$ if $(n,f(0)\Delta)=1$,
and it is declared to be a probable prime by the following algorithm.
(Such an integer will be called a {\bf Frobenius probable prime} with
respect to $f(x)$.)
All computations are done in $(\Bbb{Z}/n\Bbb{Z})[x]$.

{\bf Factorization Step} Let $f_0(x)=f(x) \bmod n$.
For $1\le i\le d$, let $F_i(x)=\gcmd(x^{n^i}-x,f_{i-1}(x))$ and
$f_i(x)=f_{i-1}(x)/F_i(x)$.
If any of the $\gcmd$s fail to exist, declare $n$ to be composite and stop.
If $f_d(x)\not=1$, declare $n$ to be composite and stop.

{\bf Frobenius Step}  For $2\le i\le d$, compute $F_i(x^n) \bmod F_i(x)$.  If it is nonzero for some $i$, 
declare $n$ to be composite and stop.

{\bf Jacobi Step} Let $S=\sum\Sb{2|i}\endSb\deg(F_i(x))/i$.

If $(-1)^S\not=\left(\Delta\over n\right)$,
declare $n$ to be composite and stop.

\vskip.05in

If $n$ is not declared to be composite by one of these three steps,
 declare $n$ to be a Frobenius probable prime and stop.
\endproclaim

The Factorization Step produces a ``distinct degree'' factorization
when $n$ is prime.  It may be of some interest to apply algorithms that
factor the polynomials completely, thus developing definitions for Berlekamp
and Cantor-Zassenhaus probable primes.  
The Cantor-Zassenhaus algorithm shares ideas with the ``strong''
Frobenius probable prime test of Section 5.  Berlekamp's algorithm has
two forms; one deterministic and one probabilistic.  The deterministic
version has running time proportional to $n$, so it is too slow to be
used in primality testing.  The probabilistic version is fast, but
since it is significantly more complicated than most existing probable
prime tests, we omit consideration of it here.

\proclaim{Corollary 3.6}
Every odd prime $p$ is a Frobenius probable prime with respect to any
monic polynomial $f(x)$ such that $p$ does not divide $f(0)\Delta$.
\endproclaim
\demo{Proof}
Immediate from Theorem 3.1.
\enddemo

\proclaim{Definition}
Let $f(x)\in\Bbb{Z}[x]$.
A {\bf Frobenius pseudoprime} with respect to a monic polynomial
$f(x)$ is a composite which is a Frobenius probable
prime with respect to $f(x)$.
\endproclaim

\head{\S 4 The Relation of Frobenius Pseudoprimes to Other
Pseudoprimes}
\endhead
\proclaim{Theorem 4.1}
An odd integer $n$ is a pseudoprime to the base $a$ if and only if it is a
Frobenius pseudoprime with respect to the polynomial $f(x)=x-a$.
\endproclaim
\demo{Proof}
First, assume $n$ is an pseudoprime to the base $a$. Thus, $a^{n-1}\equiv 1 \bmod n$.  Since $\Delta=1$,
$(n,f(0)\Delta)=(n,a)=1$.
Because
$(x-a)|(x^n-a^n)$ and $a^n-a=0$ in $(\Bbb{Z}/n\Bbb{Z})[x]$,
we have $(x-a)|((x^n-a^n)+(a^n-a)+(a-x))$ or
$(x-a)|(x^n-x)$.  Therefore, $F_1(x)=x-a$, and $f_1(x)=1$, so $n$
passes the Factorization Step.  Since $d=1$, The Frobenius Step is vacuous.
Note that $S=0$ and 
$\left(\Delta\over n\right)=\left(1\over n\right)=1$,
so $n$ passes the Jacobi Step.  Therefore $n$ is a
Frobenius pseudoprime.

Now, assume $n$ is a Frobenius pseudoprime with respect to $x-a$.  In
order to have $f_1(x)=1$, we must have $F_1(x)=x-a$.
 So by the Factorization Step,
 $(x-a)|(x^n-x)$.  Since $(x-a)|(x^n-a^n)$, $(x-a)$ divides $a^n-a$.  Since
the latter is a constant and $x-a$ is monic, $a^n-a$
 must be $0$ in $(\Bbb{Z}/n\Bbb{Z})[x]$.  So $a^n\equiv a \bmod n$.
Since $(n,f(0)\Delta)=(n,a)=1$, we have $a^{n-1}\equiv 1 \bmod n$.
Thus $n$ is a pseudoprime to the base~$a$.
\enddemo

In fact, a more general result can be proven: a Frobenius
pseudoprime with respect to $f(x)$ is a pseudoprime to the base
$f(0)$.
The idea for the proof can be found in \cite{\gurak}, Corollary 1.
We first need to prove a lemma about polynomials mod $p^k$.

\proclaim{Lemma 4.2}
Let $g(y)$ be a polynomial in $\Bbb{Z}[y]$, irreducible mod $p$.  Let
$f(x)$ be a polynomial in $\Bbb{Z}[x]$ with $p\nmid\disc(f)$.
If $f(x)$
 has $d$ roots in $(\Bbb{Z}[y]/(p,g(y)))[x]$, then it has $d$ roots
in $(\Bbb{Z}[y]/(p^k,g(y)))[x]$ for $k$ a positive integer.
\endproclaim
\demo{Proof}
By Hensel's Lemma, a root mod $p$ lifts to exactly one root mod
$p^k$, since the discriminant of $f$ is non-zero mod $p$.
Hensel's Lemma applies to any finite field of characteristic $p$.
\enddemo

\proclaim{Theorem 4.3}
Let $f(x)$ be a monic, squarefree polynomial in $\Bbb{Z}[x]$.  If an
odd integer $n$
is a Frobenius pseudoprime with respect to $f(x)$, then it is a
pseudoprime to the base $f(0)$.  
\endproclaim
\demo{Proof}
It suffices to prove $f(0)^n\equiv f(0) \bmod p^k$ for every prime
power $p^k|n$.

Let $d$ be the degree of $f(x)$.  There exists an extension field of
$\Bbb{F}_p$,
$\Bbb{F}_{p}[y]/(g(y))$, in which $f(x)$ splits completely.  The $d$ roots
must be distinct, since $n$ is coprime to the discriminant of $f(x)$.
Thus there are $d$ distinct roots of $f(x)$ in
$\Bbb{Z}[y]/(p^k,g(y))$, by Lemma 4.2.

Call the roots $y_1,y_2,\dots,y_d$.  Consider the map
$y_i\mapsto y_i^n$.  By the Frobenius Step, this map sends each root to
another root.
By the Factorization Step, $\left(y_i^n\right)^{n^{d-1}}=y_i$,
so the map is invertible.  Therefore, it permutes the roots.
Thus $\prod_{i=1}^d y_i=\prod_{i=1}^d y_i^n$.  But $\prod_{i=1}^d
y_i\equiv (-1)^df(0) \bmod p$, so $(-1)^df(0)\equiv ((-1)^df(0))^n
\bmod p^k$.  Simplifying, we see that $f(0)^n\equiv f(0) \bmod p^k$ for
each $p^k|n$.  Thus $n$ is a pseudoprime to the base $f(0)$.
\enddemo

Theorem 4.3 can be used, in conjunction with results about the
distribution of pseudoprimes base $a$, to give an upper bound on the
number of Frobenius pseudoprimes with respect to a given polynomial
$f(x)$.

\proclaim{Corollary 4.4}
Let $f(x)\in\Bbb{Z}[x]$ be a monic polynomial with nonzero
discriminant.  If $|f(0)|\not=1$, then the number of Frobenius pseudoprimes 
with respect to
$f(x)$ up to $y$ is less than $y^{1-\log\log\log y/2\log\log y}$,
for $y$ sufficiently large, where ``sufficiently large'' depends only
on $|f(0)|$.
\endproclaim
\demo{Proof}
Immediate from Theorem 4.3 and \cite{\pom}.
\enddemo

When $|f(0)|=1$, it is possible that all integers are Frobenius
pseudoprimes with respect to $f(x)$, such as if $f(x)=x-1$.  
In fact, we conjecture that for every monic,
squarefree polynomial $f(x)$, not
the product of cyclotomic polynomials, the bound of Corollary 4.4 holds.
For quadratic polynomials, the conjecture follows from \cite{\gordpom}.

\proclaim{Theorem 4.5}
If $a,\ b$ are integers, $f(x)=(x-a)(x-b)$, and $n$ is a Frobenius pseudoprime with respect to $f(x)$,
then $n$ is a pseudoprime to both bases $a$ and $b$.
\endproclaim

\demo{Proof}
Since $f(x)$ factors, its discriminant must be a square, so
$\left(\Delta\over n\right)=1$.  Therefore $F_2(x)=1$ and
$F_1(x)=f(x)$ by the Jacobi Step.
  Since $f(x)|(x^n-x)$, we have $(x-a)|(x^n-x)$.
Therefore, as in the proof of Theorem 4.1, 
we conclude $n$ is a pseudoprime to the base $a$ and, similarly, base $b$.
\enddemo

\proclaim{Theorem 4.6}
If $f(x),g(x)\in\Bbb{Z}[x]$ with $(n,\disc(fg))=1$
and $n$ is a Frobenius pseudoprime with respect to $f(x)$ and $g(x)$,
then it is a Frobenius pseudoprime with respect to $f(x)g(x)$.
\endproclaim
\demo{Proof}
Let $h(x)=f(x)g(x)$.  Let $f_i(x),g_i(x),h_i(x)$ and
$F_i(x),G_i(x),H_i(x)$ be the polynomials produced in the
Factorization Steps for $f(x),g(x),h(x)$, respectively.
If a polynomial is not defined in that step (e.g., $f_{d+1}(x)$, if
$f(x)$ has degree $d$), define it to be $1$.

We will show by induction on $i$
 that $h_i(x)=f_i(x)g_i(x)$ and $H_i(x)=F_i(x)G_i(x)$.

We have that $h_0(x)=h(x)=f(x)g(x)=f_0(x)g_0(x)$.
Assume that $h_{k-1}(x)=f_{k-1}(x)g_{k-1}(x)$.  By definition,
we have $H_k(x)=\gcmd(x^{n^k}-x,h_{k-1}(x))$, should this $\gcmd$
exist.  Since $F_k(x)|f_{k-1}(x)$,
$G_k(x)|g_{k-1}(x)$, $F_k(x)G_k(x)|h_{k-1}(x)$.

Because $\gcd(\disc(fg),n)=1$, we have that
$\gcd(F_k(x),G_k(x))=1$ in $\Bbb{F}_p[x]$ for each $p|n$.
By Proposition 3.4, $\gcmd(F_k(x),G_k(x))=1$ in
$(\Bbb{Z}/n\Bbb{Z})[x]$.
Therefore $F_k(x)A_1(x)+G_k(x)A_2(x)\equiv 1\bmod n$ for some $A_1(x)$ and
$A_2(x)$ in $\Bbb{Z}[x]$.  Also, $x^{n^k}-x\equiv B_1(x)F_k(x)\equiv
B_2(x)G_k(x)$, for some $B_1(x)$ and $B_2(x)$ in $\Bbb{Z}[x]$,
 by the definitions of $F_k(x)$ and $G_k(x)$.  Thus
$F_k(x)B_1(x)\equiv G_k(x)B_2(x)\bmod n$.  Multiplying by $A_1(x)$
gives $F_k(x)A_1(x)B_1(x)\equiv G_k(x)A_1(x)B_2(x)$.  If we substitute
for
$F_k(x)A_1(x)$, we get 
$(1-G_k(x)A_2(x))B_1(x)\equiv G_k(x)A_1(x)B_2(x)$, or
$B_1(x)\equiv G_k(x)[A_2(x)B_1(x)+A_1(x)B_2(x)]$.  Hence
$G_k(x)|B_1(x)$, and
$F_k(x)G_k(x)|x^{n^k}-x$.  Thus $F_k(x)G_k(x)|H_k(x)$.

We have, by the definitions of $F_k(x)$ and $G_k(x)$, that $F_k(x)\equiv 
r_1(x)(x^{n^k}-x)+s_1(x)f_{k-1}(x)$ and $G_k(x)\equiv 
r_2(x)(x^{n^k}-x)+s_2(x)g_{k-1}(x)$, for some polynomials
$r_1(x)$, $r_2(x)$, $s_1(x)$, and $s_2(x)$.  Multiplying these two
congruences together, 
$F_k(x)G_k(x)\equiv r_3(x)(x^{n^k}-x)+s_1(x)s_2(x)h_{k-1}(x)$,
where 
$r_3(x)=r_1(x)r_2(x)(x^{n^k}-x)+r_1(x)s_2(x)g_{k-1}(x)+r_2(x)s_1(x)f_{k-1}(x)$.
Therefore, $H_k(x)=F_k(x)G_k(x)$.  Now $h_k(x)=h_{k-1}(x)/H_{k-1}(x)$, so by
the inductive hypothesis, 
$$h_k(x)=f_{k-1}(x)g_{k-1}(x)/(F_{k-1}(x)G_{k-1}(x))=
f_k(x)g_k(x).$$

Each of the $\gcmd$s in the Factorization Step exists, and
$h_{\deg(fg)}(x)=1$. Thus $n$ passes the Factorization Step.

Since $H_i(x)=F_i(x)G_i(x)$, and $n$ passes the
Frobenius Step for $f(x)$ and $g(x)$, $H_i(x)|H_i(x^n)$. 
Thus $n$ passes the Frobenius Step for $h(x)$.

Let $S_f$ and $S_g$ be the values of $S$ computed in the Jacobi Step
for $f(x)$ and $g(x)$, respectively.  Then $S=S_f+S_g$.  So it
suffices to show that $\left(\disc(fg)\over n\right)=
\left(\disc(f)\over n\right)\left(\disc(g)\over n\right)$.  To show
this equality, it suffices to show that $\disc(fg)=\disc(f)\disc(g)\ell^2$,
where $\ell\in\Bbb{Z}$.

Let $\alpha_1,\dots,\alpha_j$ be the roots of $f(x)$ and
$\alpha_{j+1},\dots,\alpha_d$ be the roots of $g(x)$, all in
$\bar\Bbb{Q}$.  Then
$\disc(fg)=\prod_{i<i'}(\alpha_i-\alpha_{i'})^2=\disc(f)\disc(g)
\prod_{i\le j<i'}(\alpha_i-\alpha_{i'})^2$.  Let $\ell=\prod_{i\le
j<{i'}}(\alpha_i-\alpha_{i'})$.  

We will show that $\ell\in\Bbb{Z}$, by showing
$\sigma(\ell)=\ell$ for all
$\sigma\in\text{Gal}(\Bbb{\bar Q}/\Bbb{Q})$.
Any such $\sigma$ must map each $\alpha_i$ with $i\le j$ to some
$\alpha_{\bar i}$ with $\bar i\le j$, and similarly for $i>j$.  Thus $\sigma$
must only rearrange terms in the product, and $\sigma(\ell)=\ell$.
Thus $\left(\disc(fg)\over n\right)=
\left(\disc(f)\over n\right)\left(\disc(g)\over n\right)$, $n$ passes
the Jacobi Step, and $n$ is a Frobenius pseudoprime with respect to
$f(x)g(x)$.
\enddemo

The converse to Theorem 4.6 is true for the product of two linear
polynomials, as Theorems 4.1 and 4.5
show.  It is not, however, true in general.  If
$f(x)=(x-1341)(x-513)(x-545)$, then $1537$ is a Frobenius pseudoprime
with respect to $f(x)$, but it is not a pseudoprime to any of the
bases $1341$, $513$, or $545$.  This example appears in
\cite{\as} and indicates the possible usefulness of the quadratic
forms test contained therein.  The examples produced in that paper, however,
all involve polynomials with relatively large discriminant compared to
the pseudoprimes.

\proclaim{Corollary 4.7}
If $n$ is a Carmichael number, $f(x)\in\Bbb{Z}[x]$ is monic, 
$f(x)$ factors into linear
factors mod $n$ and $(n,f(0)\Delta)=1$, then $n$ is a Frobenius
pseudoprime with respect to $f(x)$.
\endproclaim
\demo{Proof}
Apply Theorem 4.1 and Theorem 4.6.
\enddemo

\proclaim{Lemma 4.8}
Let $m,n$ be positive integers, and let $f(x),g(x),r(x)\in\Bbb{Z}[x]$.
If $f(r(x))\equiv 0 \bmod (n,f(x))$ and
$x^m \equiv g(x) \bmod (n,f(x))$, then
$r(x)^m \equiv g(r(x)) \bmod (n,f(x))$.
\endproclaim
\demo{Proof}
$x^m\equiv g(x)+f(x)h(x) \bmod n$, for some $h(x)\in\Bbb{Z}[x]$.
Since $x$ is an
indeterminate, $r(x)^m\equiv g(r(x))+f(r(x))h(r(x)) \bmod n$. Because
$f(r(x))\equiv 0 \bmod (n,f(x))$, $r(x)^m\equiv g(r(x)) \bmod (n,f(x))$.
\enddemo

\proclaim{Theorem 4.9}
If $f(x)=x^2-Px+Q\in\Bbb{Z}[x]$, and $n$ is a Frobenius 
pseudoprime with respect to $f(x)$, then $n$ 
is a Lucas pseudoprime with parameters $(P,Q)$.
\endproclaim
\demo{Proof}
Note that $S=0$ or $1$.

If $\left(\Delta\over n\right)=1$, then we must have $S=0$, so
$x^n\equiv x \bmod (n,f(x))$.  Since $Q$ is invertible mod $n$, $x$ is
invertible mod $(n,f(x))$.  Thus $x^{n-1}\equiv 1$.  
By Lemma 4.8, $(P-x)^{n-1}\equiv 1$, since $f(P-x)\equiv 0 \bmod (n,f(x))$.

The two roots of $f(x)$ in $\Bbb{Z}[x]/(f(x))$ are $x$ and $P-x$, and
$(x-(P-x))^2\equiv P^2-4Q \bmod f(x)$.  Since $n$ is coprime to the discriminant
$P^2-4Q$, the difference of the two roots is invertible.
Thus $U_{n-1}\equiv \frac{x^{n-1}-(P-x)^{n-1}}{x-(P-x)}\equiv \frac{1-1}{2x-P}=0 \bmod (n,f(x))$.

If $\left(\Delta\over n\right)=-1$, then we must have $S=1$, so
$x^n\not\equiv x \bmod (n,f(x))$.  We cannot have $x^n\equiv x \bmod
(p^k,f(x))$ for $p^k|n$, since then $x^n-x\equiv 0 \bmod (p^k,f(x))$, and
the gcmd in the Factorization Step
would not exist.  Further, this shows that $f(x)$ is irreducible mod
$p$.

Because $p\nmid\Delta$, there are only 2 roots to $f(x) \bmod (p^k,f(x))$, by Lemma 4.2.
Since they are known to be $x$ and $P-x$, we must have
$x^n\equiv P-x \bmod (p^k,f(x))$ for each prime power $p^k|n$, by the
Frobenius Step.
Since $f(x)$ is monic, the congruence must hold mod $(n,f(x))$ by the
Chinese Remainder Theorem.
By Lemma 4.8, $(P-x)^n\equiv x \bmod (n,f(x))$, so $U_{n+1}\equiv\frac{x^{n+1}-
(P-x)^{n+1}}{2x-P}\equiv \frac{x(P-x)-(P-x)x}{2x-P}=0$.
\enddemo

Note that the Frobenius test is in fact stronger than the Lucas test.
For example, $323$ is the first Lucas pseudoprime with respect to the Fibonacci
sequence.  If we compute $x^{323}-x \bmod (323,x^2-x-1)$, we get $-1$.  So
$F_1(x)=1$.  If we compute $x^{323^2}-x \bmod (323,x^2-x-1)$, we get $0$.
So $F_2(x)=x^2-x-1$ and $f_2(x)=1$.
So $323$ passes the Factorization Step.  Note that it also passes the
Jacobi Step,
since $\left({5}\over 323\right)=-1$.  But it fails the Frobenius Step, because
$x^{323}\equiv x-1 \bmod (323,x^2-x-1)$, and $F_2(x-1)=-2x+2$.
The first Frobenius pseudoprime with respect to the Fibonacci
polynomial $x^2-x-1$ is $5777$.

\proclaim{Theorem 4.10}
If $f(x)=x^3-rx^2+sx-1$, then any Frobenius pseudoprime $n$ with respect
to $f(x)$ is also a Perrin pseudoprime.
In particular, if $F_1(x)=f(x)$, then $n$ has an S-signature, if
$F_3(x)=f(x)$, then $n$ has an I-signature, and if $\deg(F_1)=1$ and
$\deg(F_2)=2$, $n$ has a Q-signature.
\endproclaim
\demo{Proof}
The idea behind this proof is that relationships between $n$th powers
of the roots determine the signature,
and the necessary relationships are guaranteed to hold
because $n$ passes the Frobenius Probable Prime Test.

To this end, we use Lemma 2 of \cite{\as}.  Let $K$ be the splitting
field of $f(x)$.  Let $\alpha_1$, $\alpha_2$ and $\alpha_3$ be the three
roots of $f(x)$ in $K$.  
Lemma 2 says that $n$ has a Q-~signature if for each prime power $p^k|n$,
and for each prime ideal $\frak{p}$ of $K$ with $\frak{p}|p$, $\alpha_1^n\equiv \alpha_1$,
$\alpha_2^n\equiv\alpha_3$, and $\alpha_3^n\equiv\alpha_2 \bmod \frak{p}^k$ (or some other
permutation of the roots of order 2.)

If $\deg(F_1)=1$ and $\deg(F_2)=2$, then we must have
$f(x)\equiv F_1(x)F_2(x) \bmod p^k$.  So $f(\alpha_i)=0\equiv
F_1(\alpha_i)F_2(\alpha_i)$ for $i=1,2,3$.  Because $F_1(x)$ is linear
it has exactly one root mod $p^k$.
Therefore, one of the roots (say $\alpha_1$) is a root of $F_1(x)$ and
the other two are roots of $F_2(x)$.

For $\alpha_1$, we have $x^n\equiv x \bmod (p^k,F_1(x))$.  But we must
have $F_1(x)\equiv x-\alpha_1 \bmod p^k$, so $\alpha_1^n\equiv\alpha_1 \bmod p^k$, and hence for every prime ideal power  dividing $p^k$.

We have $x^n\not\equiv x \bmod (p^k,F_2(x))$, but $F_2(x^n)\equiv 0$.
Since there are only two roots of $F_2(x)$ mod $p^k$, 
we must have $\alpha_2^n\equiv
\alpha_3 \bmod p^k$, and similarly $\alpha_3^n\equiv\alpha_2$.

The proofs of the S and I cases are similar.
\enddemo

\proclaim{Theorem 4.11}
Let $f(x)\in\Bbb{Z}[x]$ be a monic, squarefree polynomial.
Let $\beta_1,\dots,\beta_d$
be its roots, and let $V_k=\beta_1^k+\dots+\beta_d^k$.
If $n$ is a Frobenius pseudoprime with respect to $f(x)$,
then $n$ is a pseudoprime with
respect to $V$, in the sense of \cite{\gurak}, Section 4.
\endproclaim
\demo{Proof}
The theorem follows directly from Theorem 2 of \cite{\gurak}.
\enddemo

\proclaim{Theorem 4.12}
Let $f(x)$ be a monic, squarefree polynomial.  If $n$ is a Frobenius
pseudoprime with respect to $f(x)$, then $n$ is a pseudoprime in the
sense of Szekeres.
\endproclaim
\demo{Proof}
It suffices to show that the map $x\mapsto\topsmash{x^n}$ permutes the roots
of $f(x)$.  This fact follows from the Frobenius Step.
\enddemo

Having presented the definition of Frobenius pseudoprime as a
generalization of other definitions of pseudoprime, we would like to
use the above theorems to produce a theorem that holds for all of these
types of pseudoprimes.

\proclaim{Conjecture 4.13}
For any monic, squarefree polynomial $f(x)\in\Bbb{Z}[x]$, 
there are infinitely many
Frobenius pseudoprimes with respect to $f(x)$.  In fact, for any
$\epsilon>0$, there exists a $T$ (depending on $f(x)$ and $\epsilon$)
such that if $t>T$, there are at least
$t^{1-\epsilon}$ Frobenius pseudoprimes less than $t$.
\endproclaim

It is straightforward to prove the first assertion for many polynomials (those
which split into linear and quadratic factors over $\Bbb{Z}$).  The
proof uses Corollary 4.7 and 
an extension of results in \cite{\agp} and
\cite{\agptwo}.  It is possible to prove this statement for all polynomials, but
the proof requires results about L-functions over number fields.  The
proof is given in \cite{\mepseudothree}.
The second assertion seems considerably more difficult to prove; 
for a discussion of impediments, see
\cite{\agp}.

\head{\S 5 Strong Frobenius Pseudoprimes}
\endhead
We can strengthen the test developed in the previous section by using
the identity $x^{n^i-1}-1=(x^s-1)\botsmash{\prod_{j=1}^r (x^{2^{j-1}s}+~1)}$
(where $n^i-1=2^rs$) to further factor $F_i(x)$.
 
\proclaim{Theorem 5.1}
Let $f(x)$, $d$, $\Delta$, $p$, and $F_i(x)$ be as in Theorem 3.1.
Let $p^i-1=2^rs$ with $s$ odd.
Let $F_{i,0}(x)=\gcd(F_i(x),x^s-1)$.
For $1\le j\le r$, let $F_{i,j}(x)=\gcd(F_i(x),x^{2^{j-1}s}+1)$.
Then $\botsmash{\prod_{j=0}^r F_{i,j}(x)=F_i(x)}$ and for each $j$, the
degree of $F_{i,j}(x)$ is divisible by $i$.
\endproclaim
\demo{Proof}
We have the identity $\topsmash{x^{p^i-1}-1=(x^s-1)\prod_{j=1}^r(x^{2^{j-1}s}+1)}$.
The result follows since the factors in the product are pairwise
coprime, and since $f(0)\ne 0$.
\enddemo
\proclaim{Definition}
Let $f(x)\in\Bbb{Z}[x]$ be a monic polynomial of degree $d$ with discriminant
$\Delta$.  An odd integer $n$ with $(n,f(0)\Delta)=1$ is said to pass
the {\bf strong Frobenius probable prime test} with respect to $f(x)$
if it is a Frobenius probable prime and is declared to be a probable prime
by the following additional step. (Such an integer will be called a
{\bf strong Frobenius probable prime} with respect to $f(x)$.)

{\bf Square Root Step} For each $1\le i\le d$, let $n^i-1=2^rs$ with $r$
odd.  
Let ${F}_{i,0}(x)=\gcmd(F_i(x),x^s-1)$.
Let ${F}_{i,j}(x)=\gcmd({F}_i(x), x^{2^{j-1}s}+1)$.  
Then if 
${F}_i(x)\neq\prod_{j=0}^r F_{i,j}(x)$,
if for some $j$, the degree of $F_{i,j}(x)$ is not a
multiple of $i$, or if one of the $\gcmd$s fails to exist,
declare $n$ to be composite and terminate.

If $n$ is not declared to be composite by the Frobenius probable prime test
or the Square Root Step, declare $n$ to be a strong Frobenius probable prime.
\endproclaim

\proclaim{Corollary 5.2}
Every odd prime $p$
is a strong Frobenius probable prime with respect to any monic
polynomial $f(x)$ such that $p$ does not divide $f(0)\Delta$.
\endproclaim

\proclaim{Definition}
A {\bf strong Frobenius pseudoprime} with respect to a monic polynomial
$f(x)\in\Bbb{Z}[x]$ is a composite strong Frobenius probable prime with
respect to $f(x)$.
\endproclaim

Clearly every strong Frobenius pseudoprime with respect to $f(x)$
is a Frobenius pseudoprime with respect to $f(x)$.

\proclaim{Theorem 5.3}
A number $n$ with $(n,2a)=1$ is a strong Frobenius pseudoprime
with respect to $x-a$ if and only if $n$
 is a strong pseudoprime to the base $a$.
\endproclaim
\demo{Proof}
From Theorem 3.1 it suffices to show that a pseudoprime to the base $a$ is
strong if and only if it passes the Square Root Step with respect to $x-a$.

In order to pass the Square Root Step, we need to have
$x-a|x^{2^{r-j}s}+1$ for some $1\le j\le r$ or $x-a|x^s-1$.  The first
statement is equivalent to $a^{2^{r-j}s}\equiv -1 \bmod n$ and the
second is equivalent to $a^s\equiv 1 \bmod n$.  These are exactly the
conditions for strong pseudoprimality. So $n$ passes the Square Root
Step if and only if it is a strong pseudoprime to the base $a$.
\enddemo

\proclaim{Corollary 5.4}
Every strong Frobenius pseudoprime
with respect to $x-a$  is an Euler pseudoprime to the base $a$.
\endproclaim

The situation with strong Lucas pseudoprimes is a bit more complicated, as
the polynomial needs to be changed.

\proclaim{Theorem 5.5}
Let $f(x)=x^2-Px+Q$.
Let $n$ be a integer with $(n,2\Delta Q)=1$.  Let $Q'$ be an integer
with $Q'\equiv Q^{-1} \bmod n$.
If $n$ is a strong Frobenius pseudoprime with respect to $X^2+(2+b^2c')X+1$,
then $n$
is a strong Lucas pseudoprime with parameters $(P,Q)$.
\endproclaim
\demo{Proof}
Let $U_k=U_k(P,Q)$ and $V_k=V_k(P,Q)$.
Note that $U_k\equiv 0 \bmod n$ if and only if $x^k-(P-x)^k\equiv 0 \bmod 
(n,x^2-Px+Q)$ if and only if
$\left(\frac{P-x}x\right)^k \equiv 1$.  Similarly, $V_k \equiv 0$ if
and only if
$\left(\frac{P-x}x\right)^k\equiv -1$.  Let $X=-Q'Px+Q'P^2-1$.  Then
$X\equiv (P-x)/x \bmod (n,x^2-Px+Q)$ and
$X^2+(2-P^2Q')X+1=(PQ')^2(x^2-Px+Q)$.  So, by a change of variables,
we see that $U_k\equiv 0 \bmod n$ if and only if $X^k\equiv 1 \bmod
\looseness=1
(n,X^2+(2-P^2Q')X+1)$.  The same statement holds for $V_k$, with $1$
replaced by $-1$.  So $U_k\equiv 0 \bmod n$ if and only if $X^2+(2-P^2Q')X+1$
divides $X^k-1$.  Using this statement, the fact that $n$ is a strong Lucas pseudoprime with parameters $(P,Q)$
follows immediately from the Square Root Step.
\enddemo

If we insist on keeping the same polynomial, a weaker result can be
proven.

\proclaim{Theorem 5.6}
Every strong
Frobenius pseudoprime $n$ with respect to $f(x)=x^2-Px+Q$ such that
$\left(P^2-4Q\over n\right)=-1$ is a strong Lucas
pseudoprime with parameters $(P,Q)$.
\endproclaim
\demo{Proof}
Let $U_k=U_k(P,Q)$ and $V_k=V_k(P,Q)$.
Write $n+1=2^RS$ and $n^2-1=2^rs$ with $s$ and $S$ odd.  Note that
$2^rs=(n-1)2^RS$, so $R<r$ and $S|s$.

Observe that the only ways to pass the Square Root Step are if
$f(x)|x^s-1$
or $f(x)|x^{2^{r-j}s}+1$ for some $j$ such that $r\ge j> 0$.

This means that either $f(x)|x^{2^{r-R}s}-1$ or $f(x)|x^{2^{r-j}s}+1$
for some $j$ such that $R\ge j>0$.

In the first case, we observe $2^{r-R}s=(n-1)S=nS-S$.  So
$x^{nS-S}\equiv 1 \bmod (n,f(x))$, or $x^{nS}\equiv x^S$.  But we know
that $x^n\equiv (P-x)$, so $(P-x)^S\equiv x^S$.  Thus $U_S \equiv
(x^S-(P-x)^S)/(x-(P-x))\equiv 0 \bmod (n,f(x))$ and thus mod $n$.

In the second case, we use the
formula $2^{r-j}s=(n-1)2^{R-j}S=n2^{R-j}S-2^{R-j}S$.
But $x^{2^{r-j}s}\equiv -1 \bmod (n,f(x))$, so $x^{2^{R-j}S}\equiv
-x^{n2^{R-j}S}$.  This gives us that $x^{2^{R-j}S}\equiv
-(P-x)^{2^{R-j}S}$.

Since $V_m\equiv x^m+(P-x)^m$, $V_{2^{R-j}S}\equiv 0 \bmod n$.
We conclude that $n$ is a strong Lucas pseudoprime.
\enddemo

Theorem 5.6 would not be true without the restriction that $\left(b^2+4c\over
n\right)=-1$.  For example, $294409$ is a strong Frobenius pseudoprime
with respect to $x^2-1185x+56437$, but it is not a strong Lucas
pseudoprime with parameters $(1185,56437)$.

\proclaim{Theorem 5.7}
If $n$ is a strong Frobenius pseudoprime with respect to $x^2-bx+1$,
then $n$ is an extra strong Lucas pseudoprime to the base $b$.
\endproclaim
\demo{Proof}
Let $U_k=U_k(b,1)$ and $V_k=V_k(b,1)$.
Assume that $\left(b^2+4\over n\right)=-1$.
Let $R,r,S,s$ be as in the proof of Theorem 5.6.
Observe that $x^{n+1}\equiv x(b-x) \equiv 1$.

If $f(x)|x^{2^{r-j}s}+1$ for some $j$ such that $r\ge j>0$, then we have that
$V_{2^{R-j}s}\equiv 0$, as in
Theorem 5.6.

If $V_{2^{r-1}s}\equiv 0$, we have that $V_{\frac{n+1}2}\equiv 0$, so
$x^{\frac{n+1}2}+(b-x)^{\frac{n+1}2}\equiv 0$.  Since $(b-x)\equiv
x^{-1}$, we deduce $x^{n+1}\equiv -1$, a contradiction.  This
establishes that $j>1$, 
as the definition of extra strong Lucas pseudoprime requires.

If $f(x)|x^s-1$, this means that $x^s\equiv 1 \bmod (n,f(x))$.  
$s=S(n-1)/2^{r-R}$.  So $\gcd(\frac sS,n+1)=1$.  Therefore
$x^S\equiv 1$, and $V_S \equiv x^S+(b-x)^S \equiv 1+1=2$,
and $U_S \equiv 0$ as above.

Similarly, if $f(x)|x^s+1$, we have $V_S\equiv -2$ and $U_S\equiv
0$.

The only remaining case is $f(x)|x^{2^{r-j}s}+1$ for some $j$ such
that  $r >j \ge R+1$.  $r>R+1$ only if $n\equiv 1 \bmod 4$.
Then $R=1$, and $j\le r-1$.  So $2^{r-j}s=(n+1)\frac{n-1}
{2^j}$, and  $x^{2^{r-j}s}\equiv x^{(n+1)\frac{n-1}{2^j}}\equiv 1$.
This contradicts the assumption that  $f(x)|x^{2^{r-j}s}+1$.

The proof for the case where the Jacobi symbol is $1$ is similar.
\enddemo

\head{\S 6 Carmichael-Frobenius Numbers}
\endhead
A Carmichael number is to be a number which is a (Fermat)
pseudoprime to every base.  With that in mind, we make the following
definition.

\proclaim{Definition}
Let $K$ be a number field and $n$ an odd composite with
$(n,\disc(K))=1$.
If, for each polynomial $f(x)\in\Bbb{Z}[x]$ with all its roots in $K$ and
$(n,f(0)\disc(f))=1$, $n$ is a Frobenius pseudoprime with respect to $f(x)$,
then $n$ is a {\bf Carmichael-Frobenius number} with respect to $K$. 
\endproclaim

Note that $n$ is a Carmichael number if and only if it is a 
Carmichael-Frobenius number with respect to $\Bbb{Q}$.  Also, if $n$
is a Carmichael-Frobenius number with respect to $K$, then it is also
a Carmichael-Frobenius number with respect to any subfield of $K$.  In
particular, a Carmichael-Frobenius number with respect to $K$ is also
a Carmichael number.

\proclaim{Proposition 6.1}
Let $n$ be a Carmichael number, and let $K$ be a number field with
$(n,\disc(K))=1$. If every prime
$p|n$ splits completely in $K$, then $n$ is a Carmichael-Frobenius number
with respect to $K$.
\endproclaim
\demo{Proof}
Let $f(x)\in\Bbb{Z}[x]$ be a polynomial with all of its roots in $K$ such
that
$\gcd(n,f(0)\disc(f))=1$.
For each $p|n$, $f(x)$ must split into linear factors mod $p$, since
$p$ splits completely in $K$.  Since $n$ is a Carmichael number, it is
squarefree, so $f(x)$ splits into linear factors mod $n$.  The
Proposition follows from Corollary 4.7.
\enddemo

These Carmichael-Frobenius numbers have $F_1(x)=f(x)$ in the
Factorization Step for each $f(x)$ with all of its roots in $K$.  In 
\cite{\mepseudothree},
we will show that there are infinitely many of them for each number
field $K$.  Other types of
Carmichael-Frobenius numbers are harder to come by.  The methods of 
\cite{\dopo} can be used to give heuristics suggesting that there are
infinitely many Carmichael-Frobenius numbers with respect to $K$
with $F_2(x)=f(x)$ for each irreducible
$f(x)$ with all of its roots in $K$.  We also have the following Proposition,
which is similar to Proposition 6 of \cite{\gurak}.

\proclaim{Proposition 6.2}
Let $f(x)\in\Bbb{Z}[x]$ be a monic, irreducible
polynomial of degree $k$ with splitting field $K$.
Let $n$ be a Carmichael-Frobenius number with respect to $K$.  If $F_k(x)=
f(x)$
in the Factorization Step of the Frobenius Probable Prime Test with
respect to $f(x)$, then $n$ has at least $k+2$ prime factors.
\endproclaim
\demo{Proof}
Let $p$ be a prime factor of $n$, and let $f_p(x)$ be an irreducible
factor of $f(x)$ of maximal degree in $\Bbb{F}_p[x]$.  
Let $A_p=\Bbb{F}_p[x]/(f_p(x))$.  We have that 
$A_p =\Bbb{F}_{p^r}$ for some $r\ge 1$.

We will show that $r=k$.  Since $x^n$ is a root of $f(x)$ in 
$\Bbb{Z}[x]/(n,f(x))$, it is a root in $A_p$, 
and we must have $x^n=x^{p^t}$ in $A_p$,   for some $t>0$.
We thus have $x^{n^r}\equiv x^{p^{tr}}\equiv x$ in $A_p$.
Thus $f_p(x)|\gcd(f(x),x^{n^r}-x)$ in
$\Bbb{F}_p[x]$.  Since all $\gcmd$s were computable,
$f_p(x)|F_{r'}(x)$, for some $r'\le r$.  But since $F_k(x)=f(x)$,
we must have $r'=k$, and thus
$r=k$.

Let $\alpha$ be a root of $f(x)$ in $K$.  Then for some
$g_p(x)\in\Bbb{Z}[x]$, $g_p(\alpha)$ has order $p^k-1$ in $A_p^*$.
By the Chinese Remainder Theorem, there is a monic 
polynomial $g(x)\in\Bbb{Z}[x]$ such that $g(x)\equiv g_p(x) \bmod p$
 for each $p|n$.  
Let $h(x)$ be the minimal polynomial of $g(\alpha)$ over $\Bbb{Q}$.
Then $h(x)$ has all of its roots in $K$.  Since $h(x)$, considered mod $p$,
is the minimal polynomial for $g_p(\alpha)$, we have $p\nmid h(0)\disc(h)$
for each $p|n$,
and thus $\gcd(n,h(0)\disc(h))=1$.

Thus $h(x^n)=0$ in $A_p$.  But the roots of $h(x)$ in $A_p$
are $x^p, x^{p^2},\dots,
x^{p^{k-1}}, x^{p^k}$.
Then $n\equiv p^t\bmod (p^k-1)$, for some $1\le t \le r$.
This congruence gives $p^k-1|n-p^t$, for some $1\le t\le k$.  Therefore,
$p^k-1|\frac np-p^{t-1}$.  Since $n$ is a Carmichael number, it is not
a prime power, and $n>p^t$, which implies $\frac np-p^{t-1}>0$. 
So $p^k-1\le \frac np-p^{t-1}\le
\frac n p-1$;  thus $p^k\le \frac np$.  Since $n\neq p^{k+1}$, $p^k<\frac np$.
Thus we have for all $p|n$, $p^k<\frac np$, or $p^{k+1}<n$.  If $n$ has
less than or equal to $k+1$ prime factors, we have $n^{k+1}<n^{k+1}$, by
taking the product over all prime factors.  The contradiction gives
the proposition.
\enddemo

\head{\S 7 Implementation Issues}
\endhead
Performing the Frobenius test as stated
 on quadratic polynomials would seem to require computing $x^{n^2}$.  
As the theorem below shows, there is an equivalent version of the test
that merely requires computing $x^n$.

\proclaim{Theorem 7.1}
Let $f(x)=x^2-bx-c$.  Let $\Delta=b^2+4c$.   Let $n$ be an integer with
$(n,2f(0)\Delta)=1$. If $\left(\Delta\over n\right)=1$ and
$x^n\equiv x \bmod (n, f(x))$,
then $n$ is a Frobenius probable prime
with respect to $f(x)$.  If $\left(\Delta\over n\right)=-1$ and
$x^n\equiv b-x \bmod (n,f(x))$, then $n$ is a Frobenius probable prime
with respect to $f(x)$.  
\endproclaim
\demo{Proof}
If $\left(\Delta\over n\right)=1$, then the fact that $f(x)$ divides
$x^n-x$
verifies both the Factorization Step and the Jacobi Step.  The
Frobenius Step is trivial.

Suppose $\left(\Delta\over n\right)=-1$, and $x^n\equiv
b-x \bmod (n,f(x))$.  By Lemma 4.8, $(b-x)^n\equiv x \bmod (n,f(x))$,
and so 
$x^{n^2}\equiv x \bmod (n,f(x)))$.
Note that $4f(b/2)=-\Delta$ is coprime to $n$, so $x^n-x\equiv b-2x$
has gcmd $1$ with $f(x)$ in $(\Bbb{Z}/n\Bbb{Z})[x]$.
Thus $F_1(x)=1$ and  $F_2(x)=f(x)$, so $n$ passes the
Factorization and Jacobi steps.
Since $f(b-x)\equiv f(x)\equiv 0 \bmod (n,f(x))$, it passes the Frobenius Step.
\enddemo

We will leave a proof of the running time and an description of how to
speed the
strong test to \cite{\mepseudotwo}.

Note that Lemma 4.8 can also be used to speed up the test with any
degree polynomial.

Also, when computing $\gcmd(x^n-x,f(x))$, the first step should be to
compute $x^n \bmod f(x)$.  Then the Euclidean algorithm can be applied
to two polynomials whose degree is at most that of $f(x)$.

Although the Square Root Step is listed as a separate step, in
practice it would be integrated into the Factorization Step.  
A description of how to do this in the quadratic case is given in
\cite{\mepseudotwo}.

\head{\S 8 A Challenge}
\endhead
Pomerance, Selfridge and Wagstaff offer \$620
for a number $2$ or $3 \bmod 5$ that is a pseudoprime to the base $2$  and 
also a Lucas pseudoprime with respect to the Fibonacci sequence
or a proof that none exists \cite{\psw, \guy}.

In this spirit, I have offered \$6.20 for a Frobenius pseudoprime with
respect to $x^2+5x+5$ that is congruent to $2$ or $3 \bmod 5$.
This polynomial is used instead of the Fibonacci polynomial because
$x^{2(p+1)}\equiv 1 \bmod (p,x^2-x-1)$, if $p$ is $2$ or $3 \bmod
5$.  With $x^2+5x+5$, there is no similar guarantee $x$ will have small
order mod $p$.

The lower monetary figure is a reflection of my financial status at
the time of the offer, not
of any lower confidence level.  Heuristics \cite{\dopo} suggest that an
example should exist for the PSW test, and these heuristics can be
modified to suggest
that it should also be possible to find one for the above Frobenius test.  I
believe that the two problems are equally challenging.  A
justification for my belief is that an $n$ that passes my challenge
must be a pseudoprime to the base $5$ (by Theorem 4.3) as well as a Lucas
pseudoprime with parameters $(-5,-5)$ (by Theorem 4.9).


\Refs
\ref\no\adams
\by W. W. Adams
\paper Characterizing pseudoprimes for third-order linear recurrence
sequences
\jour Math. Comp.
\vol 48
\yr 1987
\pages 1--15
\endref

\ref\no\as
\by W. W. Adams and D. Shanks
\paper Strong primality tests that are not sufficient
\jour Math. Comp.
\vol 39
\yr 1982
\pages 255--300
\endref

\ref\no\agp
\by W. R. Alford, Andrew Granville and Carl Pomerance
\paper There are infinitely many Carmichael numbers
\jour Annals of Mathematics
\yr 1994
\vol 140
\pages 703--722
\endref

\ref\no\agptwo
\by W. R. Alford, Andrew Granville, and Carl Pomerance
\paper On the difficulty of finding reliable witnesses
\inbook Algorithmic Number Theory
\bookinfo Lecture Notes in Comput. Sci.
\publ Springer-Verlag
\publaddr New York
\eds L. M. Adleman and M.-D. Huang
\yr 1994
\pages 1--16
\endref

\ref\no\arno
\by S. Arno
\paper A note on Perrin pseudoprimes
\jour Math. Comp.
\yr 1991
\vol 56
\pages 371--376
\endref

\ref\no\atkin
\by A. O. L. Atkin
\paper Intelligent Primality Test Offer
\inbook Computational Perspectives on Number Theory
\bookinfo Proceedings of a Conference in Honor of A. O. L. Atkin
\eds D. A. Buell and J. T. Teitelbaum
\publ International Press
\yr 1998
\pages 1--11
\endref

\ref\no\bw
\by R. Baillie and S. S. Wagstaff, Jr.
\paper Lucas pseudoprimes
\jour Math. Comp.
\vol 35
\yr 1980
\pages 1391-1417
\endref

\ref\no\gord
\by D. M. Gordon
\paper Pseudoprimes on elliptic curves
\inbook Th\'eorie des nombres
\publ de Gruyter
\publaddr Berlin
\eds J. M. DeKoninck and C. Levesque
\yr 1989
\pages 290--305
\endref

\ref\no\gordpom
\by D. M. Gordon and C. Pomerance
\paper The distribution of Lucas and elliptic pseudoprimes
\jour Math. Comp.
\vol 57
\yr 1991
\pages 825--838
\endref

\ref\no\mediss
\by J. Grantham
\paper Frobenius Pseudoprimes
\paperinfo dissertation
\publ University of Georgia
\yr 1997
\endref

\ref\no\mepseudotwo
\by J. Grantham
\paper A Probable Prime Test With High Confidence
\vol 72
\jour J. Number Theory
\pages 32--47
\yr 1998
\endref

\ref\no\mepseudothree
\by J. Grantham
\paper There Are Infinitely Many Perrin Pseudoprimes
\endref

\ref\no\gurak
\by S. Gurak
\paper Pseudoprimes for higher-order linear recurrence sequences
\jour Math. Comp.
\yr 1990
\vol 55
\pages 783--813
\endref


\ref\no\guy
\by R. K. Guy
\book Unsolved Problems in Number Theory
\bookinfo Second Edition
\yr 1994
\publ Springer-Verlag
\publaddr New York
\page 28
\endref

\ref\no\jac
\by N. Jacobson
\book Basic Algebra I
\bookinfo Second Edition
\yr 1985
\publ W.H. Freeman
\publaddr New York
\page 258
\endref


\ref\no\ksw
\by G. C. Kurtz, D. Shanks, and H. C. Williams
\paper Fast primality tests for numbers less than $50\cdot 10\sp 9$
\jour Math. Comp.
\yr 1986
\vol 46
\pages 691--701
\endref

\ref\no\len
\by H. W. Lenstra, Jr.
\paper Primality testing
\inbook Computational Methods in Number Theory
\bookinfo Part I
\eds H. W. Lenstra, Jr. and R. Tijdeman
\publ Math. Centre Tract
\vol 154
\publaddr Amsterdam
\yr 1982
\pages 55--77
\endref

\ref\no\mojones
\by Z. Mo and J. P. Jones
\paper A new primality test using Lucas sequences
\paperinfo preprint
\endref

\ref\no\mon
\by L. Monier
\paper Evaluation and comparison of two efficient probabilistic
primality testing algorithms
\jour Theoretical Computer Science
\vol 12
\yr 1980
\pages 97-108
\endref

\ref\no\dopo
\by C. Pomerance
\paper Are there counter-examples to the Baillie -- PSW primality
test?
\inbook Dopo Le Parole aangeboden aan Dr. A. K. Lenstra
\yr 1984
\publaddr Amsterdam
\eds H. W. Lenstra, Jr., J. K. Lenstra and P. Van Emde Boas
\endref

\ref\no\pom
\by C. Pomerance
\paper On the distribution of pseudoprimes
\jour Math. Comp.
\vol 37
\yr 1981
\pages 587--593
\endref

\ref\no\psw
\by C. Pomerance, J. L. Selfridge and S. S. Wagstaff, Jr.
\paper The pseudoprimes to $25\cdot 10^9$
\jour Math. Comp.
\vol 35
\yr 1980
\pages 1003-1026
\endref

\ref\no\rabin
\by M. O. Rabin
\paper Probabilistic algorithm for testing primality
\jour J. Number Theory
\vol 12
\yr 1980
\pages 128--138
\endref

\ref\no\robin
\by R. M. Robinson
\paper The converse of Fermat's theorem
\jour Amer. Math. Monthly
\vol 64
\yr 1957
\pages 703--710
\endref

\ref\no\rot
\by A. Rotkiewicz
\paper On the pseudoprimes of the form $ax+b$ with respect to the
sequence of Lehmer
\jour  Bull. Acad. Polon. Sci. S\'er. Sci. Math. Astronom. Phys.
\vol 20
\yr 1972
\pages 349--354
\endref

\ref\no\rots
\by A. Rotkiewicz
\paper On Euler Lehmer pseudoprimes and Strong Lehmer pseudoprimes
with parameters $L,Q$ in arithmetic progressions
\jour Math. Comp.
\vol 39
\yr 1982
\pages 239--247
\endref



\ref\no\sze
\by G. Szekeres
\paper Higher order pseudoprimes in primality testing
\inbook Combinatorics, Paul Erd\H os is eighty
\bookinfo Bolyai Soc. Math. Stud.
\publ J\'anos Bolyai Math Soc.
\publaddr Budapest
\yr 1996
\vol 2
\pages 451--458
\endref
\endRefs
\enddocument